\documentstyle[11pt]{article}
\newcommand{\beq}[1]{\begin{equation}\label{#1}}
\newcommand{\eq}{\end{equation}}
\newcommand{\beqn}[1]{\begin{eqnarray}\label{#1}}
\newcommand{\eqn}{\end{eqnarray}}

\newcommand{\mbesj}{I_\nu^{(j)}((1-q^2)z;q^2)}
\newcommand{\mbes}{I_\nu^{(1)}((1-q^2)z;q^2)}
\newcommand{\mmbes}{I_{-\nu}^{(1)}((1-q^2)z;q^2)}
\newcommand{\mmbesj}{I_{-\nu}^{(j)}((1-q^2)z;q^2)}
\newcommand{\mbess}{I_\nu^{(2)}((1-q^2)z;q^2)}
\newcommand{\makd}{K_\nu^{(1)}((1-q^2)z;q^2)}
\newcommand{\makdo}{K_\nu^{(2)}((1-q^2)z;q^2)}
\newcommand{\makdj}{K_\nu^{(j)}((1-q^2)z;q^2)}
\newcommand{\basgyp}{\phantom1_2\Phi_1
(q^{\nu+1/2},q^{-\nu+1/2};-q;q,\frac{2q}{(1-q^2)z})}
\newcommand{\mbasgyp}{\phantom1_2\Phi_1
(q^{\nu+1/2},q^{-\nu+1/2};-q;q,-\frac{2q}{(1-q^2)z})}
\newcommand{\qexp}{e_q(\frac{1-q^2}{2}z)}
\newcommand{\qqexp}{e_{q^2}(\frac{(1-q^2)^2}{4}z^2)}
\newcommand{\mqexp}{e_q(-\frac{1-q^2}{2}z)}
\newtheorem{predl}{Proposition}[section]
\newtheorem{defi}{Definition}[section]
\newtheorem{rem}{Remark}[section]
\newtheorem{cor}{Corollary}[section]
\newtheorem{lem}{Lemma}[section]
\begin{document}

\begin{flushright}
q-alg/9509xxxxx \\
 ITEP-TH-6/95\\
\end{flushright}
\vspace{10mm}
\begin{center}
{\Large \bf The Modified $q$-Bessel Functions and the $q$-Bessel-Macdonald
Functions}\\ \vspace{5mm}
M.A.Olshanetsky \footnote{Supported by ISF-MIF-300 , RFFI-94-02-14365
 and INTAS-93-0166 grants} \\ ITEP,117259,
Moscow, Russia\\e-mail olshanez@vxdesy.desy.de\\
\vspace{5mm}
V.-B.K.Rogov \footnote{Supported by ISF grant JC 7100}\\ MIIT,101475,
 Moscow, Russia\\
e-mail m10106@sucemi.bitnet\\
September  1995\\
\end{center}

\begin{abstract}
We define a $q$-analog of the modified Bessel and Bessel-Macdonald
functions.
As for the q-Bessel functions of Jackson there is a couple of functions of
 the both kind. They are arisen in the Harmonic
analysis on quantum symmetric spaces similarly to their classical counterpart.
 Their definition is based on the power expansions. We derive the recurrence
relations, difference equations, q-Wronskians, and an analog of asymptotic
expansions which turns out is exact in some domain if $q\neq 1$. Some relations
for the basic hypergeometric function which follow from this fact are
discussed.
\end{abstract}

\section{Introduction}
The q-analogues of the Bessel functions introduced ninety years ago by Jackson
\cite{Ja} are a subject of investigations in the last years \cite{V,VK}. In
these works
their properties are derived in connection with the representation theory of
quantum groups as well as their classical counterparts \cite{Vi}.

Our main interest is the q-analogues of the Bessel-Macdonald function
 $K_{\nu}$ (BMF) and modified Bessel functions $I_\nu$ (MBF).
Their presence
 in the Harmonic analysis on homogenious spaces is at least threefold.
  First of all, BMF is
the essential part of the Green function for the Laplace-Beltrami operator
on noncompact symmetric spaces with nonpositive curvature \cite{D,O}. Next BMF
define irreducible representations of the isometries of the pseudoEuclidean
plane \cite{Vi}.
Finally, BMF are the simplest so-called Whittaker \cite{J,Sh,STS} functions,
 related to the symmetric spaces
of rank one. It turns out that the Whittaker functions up to a gauge
transform
coincide with the wave functions of the open Toda quantum mechanical system
\cite{STS}. In particular, for a rank one symmetric spaces they defined the
wave functions of two-dimensional gravity in the minisuperspace approximation
\cite {GM}.
In both constructions MBF play an auxiliary role .

In fact the main motivation of our work is the last approach.
In our previous work \cite{OR} we have described the Whittaker function
 for the
quantum Lobachevsky space. They emerge as some special eigenfunctions of
the quantum second Casimir operator. In this context it was possible to reduce
the problem to commutative analysis, though the original formulation is
essentially noncommutative. The further progress in the Harmonic analysis
demands to take into account
the noncommutativity of the group algebra. It concerns, in particular, with
the Poisson kernel for the Dirichlet problem
in the quantum Lobachevsky space. In the
classical case the Whittaker functions or, what is the same, BMF are the
Fourier transform of the Poisson kernel. To reproduce this construction
in the quantum case it is necessary to derive the most fundamental properties
of q-MBF.

As in the classical case we begin with definition of q-MBF as
 the power expansions.
There are two q-MBF and they are related to q-Bessel functions
 of Jackson \cite{Ja}
as the classical ones. We derive the action of difference operators,
 recurrence
relation, difference equation and q-Wronskian for them. The most of these
results can be easily derived from \cite{V,VK}, where q-Bessel functions were
investigated.
While the definition of q-MBF, based on the Jackson q-Bessel functions is
strightforward, the definition of q-BMF is a rather subtle. We choose it in a
such way that it becomes a holomorphic in the complex right half plane.
  We repeat the same program for q-BMF as for q-MBF. The most important
part for the applications is the Laurent type expansion which is only
asymptotics in the classical case and is represented by an convergent series
in some domain in the quantum situation.  Using this property  we obtain in
conclusion some relations for the basic hypergeometric function.

\section{Some preliminary relations}
We derive some needed relations for the fundamental functions $e_q(z),
E_q(z),$ and $ \Gamma_q(z)$  \cite{GR}. As by product we define
q-psi function.

{\bf1.}For an arbitrary $a$
$$
(a,q)_n=\left\{
\begin{array}{lcccl}
1 &for& n&=&0\\
(1-a)(1-aq)\cdots(1-aq^{n-1})& for& n&\ge &1,\\
\end{array}
\right.
$$
$$
(a,q)_\infty=\lim_{n\to\infty}(a,q)_n.
$$
The $q$-exponentials are determined by formulas:
\begin{equation}
e_q(z)=\sum_{n=0}^\infty\frac{z^n}{(q,q)_n}=
\frac{1}{(z,q)_\infty},\qquad |z|<1,
\label{1}
\end{equation}
\begin{equation}
E_q(z)=\sum_{n=0}^\infty\frac{q^{\frac{n(n-1)}{2}}z^n}{(q,q)_n}=
(-z,q)_\infty,\qquad.
\label{1a}
\end{equation}
They are the q-deformations of the ordinary exponent
$$
\lim_{q\to1-0}e_q((1-q)z)=\lim_{q\to1-0}E_q((1-q)z)=e^z.
$$
Obviously, we have from (\ref{1}) and (\ref{1a})
\begin{equation}
\label{1b}
e_q(z)E_q(-z)=1
\end{equation}
It follows from (\ref{1}) that the function $e_q(z)$ has the ordinary
poles in the points $z=q^{-k}, k=0,1,\ldots$.

\begin{predl}[F.I.Karpelevich]
\label{p2.1}
The $q$-exponential $e_q(z)$ can
be represented as the sum of the partial functions
\begin{equation}
\label{2}
e_q(z)=\frac{1}{(q,q)_\infty}\sum_{k=0}^\infty\frac{(-1)^kq^{k(k+1)/2}}{(q,q)_k(1-zq^k)}.
\end{equation}
\end{predl}

{\bf Proof.} Let
$$
e_q(z,n)=\frac{1}{(q,q)_n}=\sum_{k=0}^n\frac{c_{k,n}}{1-zq^k}.
$$
Then
$$
c_{k,n}=res_{z=q^{-k}}e_q(z,n)=\lim_{z\to q^{-k}}(1-zq^k)e_q(z,n)
=\frac{1}{(1-q^{-k})\cdots(1-q^{-1})(1-q)\ldots(1-q^{n-k})}=
$$
$$
=\frac{(-1)^kq^{k(k+1)/2}}{(q,q)_k(q,q)_{n-k}}.
$$
{}From this we have
$$
e_q(z,n)=\sum_{k=0}^n\frac{(-1)^kq^{k(k+1)/2}}{(q,q)_k(1-zq^k)}\frac{1}{(q,q)_{n-k}}.
$$
As $\frac{1}{(q,q)_{n-k}}<\frac{1}{(q,q)_\infty}$ we obtain
$$
e_q(z)=\lim_{n\to\infty}e_q(z,n)=\frac{1}{(q,q)_\infty}\sum_{k=0}^\infty
\frac{(-1)^kq^{k(k+1)/2}}{(q,q)_k(1-zq^k)}.\rule{5pt}{5pt}
$$
 The next Proposition follows from (\ref{1}).
\begin{predl}
\label{p2.2}
$$
\partial_qe_q((1-q)z)=e_q((1-q)z),
$$
where
\begin{equation}
\label{4}
\partial_qf(z)=\frac{f(z)-f(qz)}{(1-q)z}.
\end{equation}
\end{predl}
The following Corollaries are evident.
\begin{cor}
\label{c2.1}
\begin{equation}
\label{5}
e_q(\frac{1-q^2}{2}qz)=\left(1-\frac{1-q^2}{2}z\right)\qexp,
\end{equation}
\begin{equation}
\label{6}
e_q(\frac{1-q^2}{2}q^{-1}z)=\frac{1}{1-\frac{q^{-1}(1-q^2)}{2}z}\qexp.
\end{equation}
\end{cor}

\begin{cor}
\label{c2.2}
\begin{equation}
\label{7}
e_{q^2}(\frac{(1-q^2)^2}{4}q^2z^2)=
\left[1-\frac{(1-q^2)^2}{4}z^2\right]\qqexp.
\end{equation}
\end{cor}

\begin{rem}
\label{r2.1}
It follows from(\ref{1})
$$
\qexp\mqexp=\frac{1}{(\frac{1-q^2}{2}z,q)_\infty}
\frac{1}{(\frac{1-q^2}{2}z,q)\infty}=
\frac{1}{(\frac{(1-q^2)^2}{4}z^2,q^2)_\infty}=e_q^2(\frac{(1-q^2)^2}{4}z^2),
$$
\end{rem}

{\bf 2.} The $q$-gamma-function
\begin{equation}
\label{9}
\Gamma_q(\alpha)=\frac{q,q)_\infty}{(q^\alpha,q)\infty}(1-q)^{1-\alpha}.
\end{equation}
If $n$ is a natural number then
$$
\Gamma_q(n+1)=\frac{(q,q)_n}{(1-q)^n}.
$$

\begin{predl}
\label{p2.3}
The function $\Gamma_{q^2}(z)$ can be represanted in
the following form
\begin{equation}
\label{10}
\Gamma_{q^2}(z)=(1-q^2)^{1-z}\sum_{k=0}^\infty\frac{(-1)^kq^{k(k+1)}}
{(q^2,q^2)_k(1-q^{2k+2z})}.
\end{equation}
\end{predl}

{\bf Proof.} It follows from (\ref{9}) that
$$
\Gamma_{q^2}(z)=(q^2,q^2)_\infty(1-q^2)^{1-z}e_{q^2}(q^{2z}).
$$
Substituting in this equality (\ref{2}) we obtain (\ref{10}).
\rule{5pt}{5pt}

We denote by $\psi_q(z)$ the logarithmic derivative of
$\Gamma_q$-function
\begin{equation}
\label{11}
\psi_q(z)=\frac{\Gamma_q^\prime(z)}{\Gamma_q(z)}.
\end{equation}
{}From (\ref{9}) and (\ref{11}) we  receve immediatly

\begin{predl}
\label{p2.4}
The function $\psi_{q^2}(z)$ has the following
form
\begin{equation}
\label{12}
\psi_{q^2}(z)=-\ln(1-q^2)
+\ln{q^2}\sum_{k=1}^\infty\frac{q^{2k+2z}}{1-q^{2k+2z}}.
\end{equation}
\end{predl}

\begin{cor}
\label{c2.3}
\begin{equation}
\label{13}
\lim_{z\to-n}\frac{\psi_{q^2}(z)}{\Gamma_{q^2}(z)}=(-1)^n
\frac{q^{-n(n+1)}(q^2,q^2)_n}{(1-q^2)^{n+1}}\ln{q^2}.
\end{equation}
\end{cor}

{\bf Proof.} This result is obtained from (\ref{10}) and (\ref{12}).
\rule{5pt}{5pt}

\section{The Modify $q$-Bessel Functions}
{\bf 1. Definition.}\\
In \cite{Ja} the $q$-Bessel functions were determined as
\begin{equation}
\label{3.1}
J_\nu^{(1)}(z,q)=\frac{(q^{\nu+1},q)_\infty}{(q,q)_\infty}(z/2)^\nu
\phantom1_2\Phi_1(0,0;q^{\nu+1};q,-\frac{z^2}{4}),
\end{equation}
\begin{equation}
\label{3.2}
J_\nu^{(2)}(z,q)=\frac{(q^{\nu+1},q)_\infty}{(q,q)_\infty}(z/2)^\nu
\phantom1_0\Phi_1(-;q^{\nu+1};q,-\frac{z^2q^{\nu+1}}{4}),
\end{equation}
where $\phantom1_2\Phi_1$ and $\phantom1_0\Phi_1$ are the basic
hypergeometric functions \cite {GR}
\begin{equation}
\label{8}
\phantom1_r\Phi_s(a_1,\cdots,a_r;b_1,\cdots,b_s;q,z)=
\sum_{n=0}^\infty\frac{(a_1,q)_n\ldots(a_r,q)_n}
{(q,q)_n(b_1,q)_n\ldots(b_s,q)_n}
[(-1)^nq^{n(n-1)/2}]^{1+s-r}z^n.
\end{equation}

It allows to define q-MBF in analogy with the classical case \cite{BE}
as (\ref{3.1}) and (\ref{3.2}).

\begin{defi}
\label{d3.1}
The functions
$$
I_\nu^{(1)}(z,q)=\frac{(q^{\nu+1},q)_\infty}{(q,q)_\infty}(z/2)^\nu
\phantom1_2\Phi_1(0,0;q^{\nu+1};q,\frac{z^2}{4}),
$$
$$
I_\nu^{(2)}(z,q)=\frac{(q^{\nu+1},q)_\infty}{(q,q)_\infty}(z/2)^\nu
\phantom1_0\Phi_1(-;q^{\nu+1};q,\frac{z^2q^{\nu+1}}{4})
$$
are called the modify $q$-Bessel functions.
\end{defi}
Evidently,
$$
I_\nu^{(j)}(z,q)=e^{-\frac{i\nu\pi}{2}}J_\nu^{(j)}(e^{i\pi/2}z,q),
\qquad j=1,2.
$$
We will consider below the functions
\begin{equation}
\label{3.5}
I_\nu^{(1)}((1-q^2)z;q^2)=\sum_{k=0}^\infty\frac{(1-q^2)^k(z/2)^{\nu+2k}}
{(q^2,q^2)_k\Gamma_{q^2}(\nu+k+1)},\qquad |z|<\frac{2}{1-q^2},
\end{equation}
\begin{equation}
\label{3.6}
I_\nu^{(2)}((1-q^2)z;q^2)=
\sum_{k=0}^\infty\frac{q^{2k(\nu+k)}(1-q^2)^k(z/2)^{\nu+2k}}
{(q^2,q^2)_k\Gamma_{q^2}(\nu+k+1)}.
\end{equation}
If $|q|<1,$ the series (\ref{3.6}) converges for all $z\ne0$ absolutely.
 Therefore $\mbess$
is holomorphic function outside of $z=0$.

\begin{rem}
\label{r3.1}
$$
\lim_{q\to1-0}I_\nu^{(j)}((1-q^2)z;q^2)=I_\nu(z),\qquad j=1, 2.
$$
\end{rem}
{\bf 2. Recurrence relations and differnces.}\\
\begin{predl}
\label{p3.1}
The function $\mbes$ satisfies the following
relations
\begin{equation}
\label{3.7}
\frac{2}{(1+q)z}\partial_q{z^\nu I_{-\nu}^{(1)}((1-q^2)z;q^2)}=
z^{\nu-1}I_{-\nu+1}^{(1)}((1-q^2)z;q^2),
\end{equation}
\begin{equation}
\label{3.8}
\frac{2}{(1+q)z}\partial_qz^\nu\mbes=z^{\nu-1}I_{\nu-1}^{(1)}((1-q^2)z;q^2),
\end{equation}
where the operator $\partial_q$ is determined by (\ref{4}).
\end{predl}

{\bf Proof.}
$$
\frac{2}{(1+q)z}\partial_q{z^\nu I_{-\nu}^{(1)}((1-q^2)z;q^2)}=
\sum_{k=1}^\infty\frac{(1-q^2)^{k-1}z^{2k-2}}
{2^{-\nu+2k-1}(q^2,q^2)_{k-1}\Gamma_{q^2}(-\nu+k+1)}=
$$
$$
=\sum_{k=0}^\infty\frac{(1-q^2)^kz^{2k}}
{2^{-\nu+1+2k}(q^2,q^2)_k\Gamma_{q^2}(-\nu+1+k+1)}=
z^{\nu-1}I_{-\nu+1}((1-q^2)z;q^2).
$$
The proof of (\ref{3.8}) is the same.\rule{5pt}{5pt}

\begin{predl}
\label{p3.2}
The function $\mbes$ satisfies the following
recurrence relations
\begin{equation}
\label{3.9}
I_{\nu-1}^{(1)}((1-q^2)z;q^2)-I_{\nu+1}^{(1)}((1-q^2)z;q^2)=
\frac{2}{(1-q^2)z}(q^{-\nu}-q^\nu)I_\nu^{(1)}((1-q^2)qz;q^2),
\end{equation}
$$
I_{\nu-1}^{(1)}((1-q^2)z;q^2)+I_{\nu+1}^{(1)}((1-q^2)z;q^2)
=\frac{4}{(1-q^2)z}\mbes-
$$
\begin{equation}
\label{3.10}
-\frac{2}{(1-q^2)z}(q^{-\nu}+q^\nu)I_\nu^{(1)}((1-q^2)qz;q^2).
\end{equation}
\end{predl}

{\bf Proof.} It follows from (\ref{3.7}) and (\ref{3.8})
$$
I_{\nu-1}^{(1)}((1-q^2)z;q^2)=\frac{2}{(1-q^2)z}\left[\mbes-
q^\nu I_\nu^{(1)}((1-q^2)qz;q^2)\right],
$$
$$
I_{\nu+1}^{(1)}((1-q^2)z;q^2)=\frac{2}{(1-q^2)z}\left[\mbes-
q^{-\nu}I_\nu^{(1)}((1-q^2)qz;q^2)\right].
$$
Summing and subtracting these equalities we obtain the statement
.\rule{5pt}{5pt}

\begin{predl}
\label{p3.3}
The function $\mbess$ satisfies the following
relations
$$
\frac{2}{(1+q)z}\partial_q{z^\nu I_{-\nu}^{(2)}((1-q^2)z;q^2)}=
q^{-\nu+1}z^{\nu-1}I_{-\nu+1}^{(2)}((1-q^2)qz;q^2),
$$
$$
\frac{2}{(1+q)z}\partial_qz^\nu\mbess=
q^{-\nu+1}z^{\nu-1}I_{\nu-1}^{(2)}((1-q^2)qz;q^2).
$$
\end{predl}

\begin{predl}
\label{p3.4}
The function $\mbess$ satisfies the following
recurrence relations
$$
q^{-\nu}I_{\nu-1}^{(2)}((1-q^2)z;q^2)-q^\nu
I_{\nu+1}^{(2)}((1-q^2)z;q^2)=
\frac{2}{(1-q^2)z}(q^{-\nu}-q^\nu)I_\nu^{(2)}((1-q^2)z;q^2),
$$
\end{predl}
$$
q^{-\nu}I_{\nu-1}^{(2)}((1-q^2)z;q^2)+q^\nu
I_{\nu+1}^{(2)}((1-q^2)z;q^2)
=\frac{4}{(1-q^2)z}I_\nu^{(2)}((1-q^2)q^{-1}z;q^2)-
$$
$$
\frac{2}{(1-q^2)z}(q^{-\nu}+q^\nu)I_\nu^{(2)}((1-q^2)z;q^2).
$$

The proof of Propositions \ref{p3.3} and \ref{p3.4} are the same as
\ref{p3.1} and \ref{p3.2}.

{\bf 3. Difference equation.}\\
\begin{predl}
\label{p3.5}
The function $\mbes$ is a solution to the
difference equation
\begin{equation}
\label{3.11}
\left[1-\frac{q^{-2}(1-q^2)^2}{4}z^2\right]
f(q^{-1}z)-(q^{-\nu}+q^\nu)f(z)+f(qz)=0.
\end{equation}
\end{predl}

{\bf Proof.} Substituting (\ref{3.5}) in the left side of (\ref{3.11})
we obtain
$$
\sum_{k=0}^\infty\left[q^{-\nu-2k}-q^{-\nu}-q^\nu+q^{\nu+2k}-
\frac{q^{-\nu-2k-2}(1-q^2)^2z^2}{4}\right]\frac{(1-q^2)^k(z/2)^{\nu+2k}}
{(q^2,q^2)_k\Gamma_{q^2}(\nu+k+1)}=
$$
$$
=\sum_{k=1}^\infty\frac{q^{-\nu-2k}(1-q^{2k})(1-q^{2\nu+2k})(1-q^2)^k(z/2)^{\nu+2k}}
{(q^2,q^2)_k\Gamma_{q^2}(\nu+k+1)}-
\sum_{k=0}^\infty\frac{q^{-\nu-2k-2}(1-q^2)^{k+2}(z/2)^{\nu+2k+2}}
{(q^2,q^2)_k\Gamma_{q^2}(\nu+k+1)}.
$$
Changing $k$ to $k+1$ in the first term we obtain zero.\rule{5pt}{5pt}

\begin{cor}
\label{c3.1}
The function $I_{-\nu}^{(1)}((1-q^2)z;q^2)$ satisfies
the equation (\ref{3.11}).
\end{cor}

{\bf Proof.} Obviously the left side of (\ref{3.11}) is a even function
of
$\nu$. Then the change $\nu$ to $-\nu$ in (\ref{3.5}) transforms a
solution to (\ref{3.11}) to a solution.\rule{5pt}{5pt}

\begin{defi}
\label{d3.2}
The $q$-Wronskian of two solutions $f_\nu^1(z)$,
$f_\nu^2(z)$ to a difference second order equation is
$$
W(f_\nu^1,f_\nu^2)(z)=f_\nu^1(z)f_\nu^2(qz)-f_\nu^1(qz)f_\nu^2(z).
$$
\end{defi}

If $q$-Wronskian  nonvanishes an arbitrary solution to the
difference second order equation can be written as
$$
f_\nu(z)=C_1f_\nu^1(z)+C_2f_\nu^2(z).
$$
In this case the functions $f_\nu^1(z)$ and $f_\nu^2(z)$ is the
fundamential system of solutions.

\begin{predl}
\label{p3.6}
If $\nu$ is noninteger the functions $\mbes$ and
$I_{-\nu}^{(1)}((1-q^2)z;q^2)$ form a fundamental system of solutions
to the equation (\ref{3.11}).
\end{predl}

{\bf Proof.} Consider the $q$-Wronskian
$$
W(I_\nu^{(1)},I_{-\nu}^{(1)})(z)=
$$
\begin{equation}
\label{3.12}
\mbes I_{-\nu}^{(1)}((1-q^2)qz;q^2)-
I_\nu^{(1)}((1-q^2)qz;q^2)I_{-\nu}^{(1)}((1-q^2)z;q^2).
\end{equation}
Since $I_\nu^{(1)}$ and $I_{-\nu}^{(1)}$ are solution to (\ref{3.11})
then
$$
I_{\pm\nu}^{(1)}((1-q^2)q^2z;q^2)=
-[1-\frac{(1-q^2)^2}{4}z^2]I_{\pm\nu}^{(1)}((1-q^2)z;q^2)+
(q^{-\nu}+q^\nu)I_{\pm\nu}^{(1)}((1-q^2)qz;q^2).
$$
Thus
$$
W(I_\nu^{(1)},I_{-\nu}^{(1)})(qz)=I_\nu^{(1)}((1-q^2)qz;q^2)
I_{-\nu}^{(1)}((1-q^2)q^2z;q^2)-
$$
$$
I_\nu^{(1)}((1-q^2)q^2z;q^2)I_{-\nu}^{(1)}((1-q^2)qz;q^2)=
[1-\frac{(1-q^2)^2}{4}z^2]W(I_\nu^{(1)},I_{-\nu}^{(1)})(z).
$$
Comparing this equality with (\ref{7}) we can write
$$
W(I_\nu^{(1)},I_{-\nu}^{(1)})(z)=C_\nu\qqexp.
$$
Setting $z=0$ in (\ref{3.12}) we obtain
$$
C_\nu=\frac{q^{-\nu}(1-q^2)}{\Gamma_{q^2}(\nu)\Gamma_{q^2}(1-\nu)}.
$$
So finaly
\begin{equation}
\label{3.13}
W(I_\nu^{(1)},I_{-\nu}^{(1)})(z)=
\frac{q^{-\nu}(1-q^2)}{\Gamma_{q^2}(\nu)\Gamma_{q^2}(1-\nu)}\qqexp.
\end{equation}
Obviously this function nonvanishes.\rule{5pt}{5pt}

If $\nu=n$ is integer then from (\ref{3.5}) and (\ref{3.6})
\begin{equation}
\label{3.14}
I_{-n}^{(j)}((1-q^2)z;q^2)=I_n^{(j)}((1-q^2)z;q^2),\quad j=1, 2.
\end{equation}

\begin{predl}
\label{p3.7}
The function $\mbess$ is a solution to the
equation
\begin{equation}
\label{3.15}
f(q^{-1}z)-(q^{-\nu}+q^\nu)f(z)+\left[1-\frac{(1-q^2)^2}{4}z^2\right]f(qz)=0.
\end{equation}
\end{predl}

This Proposition is proved in the same way as Proposition \ref{p3.5}.

\begin{cor}
\label{c3.2}
 The q-MBF are related as
(\ref{3.5}) and (\ref{3.6})
\begin{equation}
\label{3.16}
\mbes=\qqexp\mbess
\end{equation}
\end{cor}

{\bf Proof.} It follows from (\ref{3.15}) and (\ref{7}) that
$$
\qqexp\mbess
$$
satisfies (\ref{3.11}). Hence if $\nu$ is noninteger
\begin{equation}
\label{3.17}
\qqexp\mbess=A\mbes+BI_{-\nu}^{(1)}((1-q^2)z;q^2).
\end{equation}
Multiplying this equality on $(z/2)^\nu$ and setting $z=0$ we obtain
$B=0$. Multiplying (\ref{3.17}) on $(z/2)^{-\nu}$ and setting $z=0$ we
obtain $A=1.$

Since (\ref{3.5}) and (\ref{3.6}) are continuous functions of $\nu$, then
(\ref{3.16}) is valid for $\nu=n$ as well.\rule{5pt}{5pt}

Multiplying the both sides of (\ref{3.16}) on
$E_{q^2}(-\frac{(1-q^2)^2}{4}z^2)$ we obtain

\begin{cor}
\label{c3.3}
$$
I_\nu^{(2)}((1-q^2)z;q^2)=E_{q^2}(-\frac{(1-q^2)^2}{4}z^2)
I_\nu^{(1)}((1-q^2)z;q^2).
$$
\end{cor}

\begin{predl}
\label{p3.8}
The function $\mbes$ is the meromorphic function
outside of $z=0$ with the ordinary poles in the points
$z=\pm\frac{2q^{-r}}{1-q^2}, r=0, 1,\ldots$
\end{predl}

{\bf Proof.} The statement of this Proposition follows from (\ref{3.16})
and Remark \ref{r2.1} as
$$
\qqexp=\frac{1}{(\frac{(1-q^2)^2}{4}z^2,q^2)_\infty}=\frac{1}
{(\frac{(1-q^2)}{2}z,q)_\infty(-\frac{(1-q^2)}{2}z,q)_\infty}=
$$
$$
=\qexp e_q(-\frac{1-q^2}{2}z).\rule{5pt}{5pt}
$$

\begin{rem}
\label{r3.2}
If $q\to1-0$ the all poles of $\mbes$
$$
z_r=\pm\frac{2q^{-r}}{1-q^2},\qquad r=0, 1,\ldots
$$
go to infinity along the real axis.
\end{rem}

\section{The Laurent Type Serieses for Modify $q$-Bessel Functions}
Unfortunately the function $\mbes$ is determined by power series
(\ref{3.5}) in domain $|z|<\frac2{1-q^2}$ only while we need a
representation of this function as series on the whole complex plane. Now we
will try to improve this situation.

\begin{predl}
\label{p4.1}
An arbitrary solution to
(\ref{3.11}) can be written in form
$$
f_\nu(z)=\frac a{\sqrt z}\qexp\basgyp+
$$
\begin{equation}
\label{4.1}
+\frac b{\sqrt z}\mqexp\mbasgyp,
\end{equation}
where $_2\Phi_1$ is determined by (\ref{8}).
\end{predl}

{\bf Proof.} Let $f_\nu(z)$ be a solution to  (\ref{3.11}).
Represent it as
\begin{equation}
\label{4.2}
f_\nu(z)=\frac a{\sqrt z}\qexp\phi_\nu^1(z)+
\frac b{\sqrt z}\mqexp\phi_\nu^2(z).
\end{equation}
We assume here that the first summand on the right side satisfyes the
equation (\ref{3.11}). Then the second summand is also a solution. Consider
the first summand. Substituting it in (\ref{3.11}) and using (\ref{5}),
(\ref{6}) we obtain the difference equation for $\phi_\nu^1(z)$
\begin{equation}
\label{4.3}
\left[1+\frac{q^{-1}(1-q^2)}{2}z\right]q\phi_\nu^1(q^{-1}z)-q^{1/2}(q^{-\nu}+q^\nu)\phi_\nu^1(z)+
\left[1-\frac{1-q^2}{2}z\right]\phi_\nu^1(qz)=0.
\end{equation}
Represent $\phi_\nu^1(z)$ as
$$
\phi_\nu^1(z)=\sum_{k=0}^\infty a_kz^{-k}.
$$
Then  $a_0$ is arbitrary, and for any $k\ge1$
$$
a_k=a_{k-1}\frac{2q(1-q^{\nu-1/2+k})(1-q^{-\nu-1/2+k})}{(1-q^2)(1-q^{2k})}.
$$
and
$$
a_k=a_0\frac{2^kq^k(q^{\nu+1/2},q)_k(q^{-\nu+1/2},q)_k}
{(1-q^2)^k(q^2,q^2)_k}.
$$
Assuming $a_0=1$ we obtain
\begin{equation}
\label{4.4}
\phi_\nu^1(z)=\sum_{k=0}^\infty
\frac{2^kq^k(q^{\nu+1/2},q)_k(q^{-\nu+1/2},q)_k}{(1-q^2)^k(q^2,q^2)_k}z^{-k}.
\end{equation}
Similarly
$$
\phi_\nu^2(z)=\sum_{k=0}^\infty(-1)^k
\frac{2^kq^k(q^{\nu+1/2},q)_k(q^{-\nu+1/2},q)_k}{(1-q^2)^k(q^2,q^2)_k}z^{-k}=
\phi_\nu^1(-z).
$$
In view of this fact we will drop the superscript.
The series (\ref{4.4}) converges absolutly for $|z|~>~\frac{2q}{1-q^2}$.

The coefficients $a$ and $b$ are determined uniquely by
$f_\nu(z)$. In fact let $z=z_0, |z_0|~>~2/(1-~q^2), f_\nu(z_0)=A$ and
$f_\nu(qz_0)=B$. Using (\ref{5}), (\ref{6}) we come to the system
for $a$ and $b$:
\begin{equation}
\label{4.5}
\left\{
\begin{array}{lcl}
\frac a{\sqrt{z_0}}\qexp\phi_\nu(z_0)+
\frac b{\sqrt{z_0}}\mqexp\phi_\nu(-z_0)&=&A\\
\frac a{\sqrt{qz_0}}(1-(1-q^2)/2z_0)\qexp\phi_\nu(qz_0)+
\frac b{\sqrt{qz_0}}(1+(1-q^2)/2z_0)\mqexp\phi_\nu(-qz_0)&=&B.\\
\end{array}
\right.
\end{equation}
Its determinant has the form
$$
W=\frac{q^{-1/2}}{z_0}\qqexp[\phi_\nu(z_0)\phi_\nu(-qz_0)-
\phi_\nu(qz_0)\phi_\nu(-z_0)+
$$
$$
+\frac{1-q^2}{2}z_0(\phi_\nu(z_0)\phi_\nu(-qz_0)+
\phi_\nu(qz_0)\phi_\nu(-z_0))].
$$
Assume for a moment that $W(z_0)=0$ for  some $z_0:|z|>\frac2{1-q^2}$.
Then
$$
\left\{
\begin{array}{rcl}
\phi_\nu(z_0)\phi_\nu(-qz_0)-\phi_\nu(qz_0)\phi_\nu(-z_0)=0\\
\phi_\nu(z_0)\phi_\nu(-qz_0)+\phi_\nu(qz_0)\phi_\nu(-z_0)=0\\
\end{array}
\right.
$$
or

\begin{equation}
\label{4.6}
\left\{
\begin{array}{rcl}
\phi_\nu(z_0)\phi_\nu(-qz_0)=0\\
\phi_\nu(qz_0)\phi_\nu(-z_0)=0\\
\end{array}
\right.
\end{equation}
It follows from (\ref{4.4}) that (\ref{4.6}) is fulfild if
$\phi_\nu(z_0)=\phi_\nu(qz_0)=0$ (or
$\phi_\nu(-z_0)=\phi_\nu(-qz_0)=0$). In this case from (\ref{4.3}) we
have $\phi_\nu(q^{-1}z_0)=0$. And hence $\phi_\nu(q^{-r}z_0)=0$ for any
$r=0, 1,\ldots$. But it contrudicts the obvious equality
$
\lim_{|z|\to\infty}\phi_\nu(z)=1.
$
Thus $W\ne0$, and $a$, $b$ are determined uniquely from (\ref{4.5}).

It is easy to see from (\ref{8}) that
$$
\phi_\nu(z)=\sum_{k=0}^\infty\frac{(q^{\nu+1/2},q)_k(q^{-\nu+1/2},q)_k}
{(q,q)_k(-q,q)_k}\left(\frac{2q}{(1-q^2)z}\right)^k=
\basgyp.
$$
Substituting this function to (\ref{4.2}) we obtain (\ref{4.1}).
\rule{5pt}{5pt}

We denote
\begin{equation}
\label{4.7}
\basgyp=\Phi_\nu(z)
\end{equation}
for brevity.

\begin{predl}
\label{p4.2}
The $q$-MBF $\mbes$ for
$z\ne0$ can be represented as
\begin{equation}
\label{4.8}
\mbes=\frac{a_\nu}{\sqrt{z}}\left[\qexp\Phi_\nu(z)+
ie^{i\nu\pi}\mqexp\Phi_\nu(-z)\right],
\end{equation}
where $\Phi_\nu(z)$ is determined by (\ref{4.7}) and
\begin{equation}
\label{4.9}
a_\nu=\sqrt{\frac{2}{1-q^2}}e_q(-1)
\frac{I_\nu^{(2)}(2;q^2)}{\Phi_\nu(\frac{2}{1-q^2})}.
\end{equation}
\end{predl}

{\bf Proof.} It follows from corollary \ref{c3.2} and Proposition
\ref{p4.1} that
$$
\qexp\mqexp\mbess=
$$
\begin{equation}
\label{4.10}
=\frac{a_\nu}{\sqrt{z}}\qexp\Phi_\nu(z)+
\frac{b_\nu}{\sqrt{z}}\mqexp\Phi_\nu(-z)
\end{equation}
in the domain $\frac{2q}{1-q^2}<|z|<\frac{2}{1-q^2}$.
The functions in the right and the left sides  are
meromorphic  in the domain $z\ne0$, and have the ordinary poles in
the points $z=\pm\frac{2q^{-r}}{1-q^2}, r=0, 1,\cdots$. Due to the
uniqueness of the analytic continuation the equality (\ref{4.10}) is
valid in the domain $z\ne0$. We require that the residues in the
poles of both sides (\ref{4.10}) are equal. Then for
$z=\frac{2q^{-r}}{1-q^2}$
\begin{equation}
\label{4.11}
e_q(-q^{-r})I_\nu^{(2)}(2q^{-r};q^2)=
a_{\nu}q^{r/2}\sqrt{\frac{1-q^2}{2}}\Phi_\nu(\frac{2q^{-r}}{1-q^2}),
\end{equation}
and for $z=-\frac{2q^{-r}}{1-q^2}$
\begin{equation}
\label{4.12}
e_q(-q^{-r})I_\nu^{(2)}(-2q^{-r};q^2)=
b_\nu\frac{q^{r/2}}{i}\sqrt{\frac{1-q^2}{2}}\Phi_\nu(-\frac{2q^{-r}}{1-q^2}).
\end{equation}
It follows from (\ref{3.6})
$I_\nu^{(2)}(-2q^{-r};q^2)=e^{i\nu\pi}I_\nu(2q^{-r};q^2)$. Hence
from (\ref{4.11}) and (\ref{4.12})
\begin{equation}
\label{4.13}
b_\nu=ie^{i\nu\pi}a_\nu.
\end{equation}
Assuming $r=0$ in (\ref{4.11}) we have (\ref{4.9}). Substituting
(\ref{4.13}) to (\ref{4.10}) we obtain the statement of the Proposition.
\rule{5pt}{5pt}

\begin{predl}
\label{p4.3}
The coefficients $a_\nu$ (\ref{4.9}) satisfy the
recurrent relation
\begin{equation}
\label{4.14}
a_{\nu+1}=a_\nu q^{-\nu-1/2},
\end{equation}
and the condition
\begin{equation}
\label{4.15}
a_\nu a_{-\nu}=\frac{q^{-\nu+1/2}}
{2\Gamma_{q^2}(\nu)\Gamma_{q^2}(1-\nu)\sin\nu\pi}.
\end{equation}
\end{predl}

{\bf Proof.} Substitute (\ref{4.8}) in (\ref{3.9}) and (\ref{3.10}).
Then using (\ref{5}) we have
\begin{equation}
\label{4.16}
a_{\nu-1}\Phi_{\nu-1}(z)-a_{\nu+1}\Phi_{\nu+1}(z)=
2a_{\nu}q^{-1/2}\frac{q^{-\nu}-q^\nu}{(1-q^2)z}(1-\frac{1-q^2}{2}z)
\Phi_\nu(qz),
\end{equation}
$$
a_{\nu-1}\Phi_{\nu-1}(z)+a_{\nu+1}\Phi_{\nu+1}(z)=
\frac{4a_\nu}{(1-q^2)z}\Phi_\nu(z)-
$$
\begin{equation}
\label{4.17}
-2a_\nu q^{-1/2}\frac{q^{-\nu}+q^\nu}
{(1-q^2)z}(1-\frac{1-q^2}{2}z)\Phi_\nu(qz).
\end{equation}
Turn $z$ to the infinity in (\ref{4.16}) and (\ref{4.17}) we come to the system
$$
\left\{
\begin{array}{rcl}
a_{\nu-1}-a_{\nu+1}=-a_{\nu}q^{-1/2}(q^{-\nu}-q^\nu)\\
a_{\nu-1}+a_{\nu+1}=a_{\nu}q^{-1/2}(q^{-\nu}+q^\nu).\\
\end{array}
\right.
$$
{}From this system we obtain the first statement of the Proposition.

Consider the $q$-Wronskian
$W(I_\nu^{(1)},I_{-\nu}^{(1)})$ (\ref{3.13}), and the representation
(\ref{4.8}). Then for $\nu\ne n$
$$
\frac{q^{-\nu}(1-q^2)}{\Gamma_{q^2}(\nu)\Gamma_{q^2}(1-\nu)}\qqexp=
q^{-1/2}a_\nu a_{-\nu}z^{-1}\qexp\mqexp\times
$$
$$
\times\{ie^{-i\nu\pi}[(1+\frac{1-q^2}{2}z)\Phi_\nu(z)\Phi_\nu(-qz)-
(1-\frac{1-q^2}{2}z)\Phi_\nu(qz)\Phi_\nu(-z)]+
$$
$$
+ie^{i\nu\pi}[(1-\frac{1-q^2}{2}z)\Phi_\nu(-z)\Phi_\nu(qz)-
(1+\frac{1-q^2}{2}z)\Phi_\nu(-qz)\Phi_\nu(z)]\}.
$$
Reducing this equality to the $q$-exponentials and turning $z$ to
infinity we obtain (\ref{4.15}).\rule{5pt}{5pt}

If $\nu=n$ from (\ref{4.14}) we have
$$
a_n=a_{n-k}q^{-k/2(2n-k)}.
$$
If $k=2n$ then $a_n=a_{-n}$.

{}From (\ref{4.15}) and \cite{GR} (1.10.16) we obtain
\begin{equation}
\label{4.18}
a_n^2=\frac{q^{-n^2+1/2}\ln{q^{-2}}}{2\pi(1-q^2)}.
\end{equation}

\begin{predl}
\label{p4.4}
The $q$-MBF $I_\nu^{(2)}((1-q^2)z;q^2)$ for $z\ne0$ can be represented
by
\begin{equation}
\label{4.19}
I_\nu^{(2)}((1-q^2)z;q^2)=\frac{a_\nu}{\sqrt{z}}
\left[E_q(\frac{1-q^2}{2}z)\Phi_\nu(z)+
ie^{i\nu\pi}E_q(-\frac{1-q^2}{2}z)\Phi_\nu(-z)\right].
\end{equation}
\end{predl}

{\bf Proof}. This statement follows from Corollary \ref{c3.3} and
(\ref{1b}).

\section{The $q$-Bessel-Macdonald Functions}

In the classical analyses the BMF are defined as
\begin{equation}
\label{5.14}
K_\nu(z)=\frac{\pi}{2\sin{\nu\pi}}\left[I_{-\nu}(z)-I_\nu(z)\right]
\end{equation}
for $\nu\ne n$, and if $\nu=n$ by the limit for $\nu\to ~n$ in
(\ref{5.14}). Here  we present  the correct "quantization" of this
definition
in a such way that other properties are also quantized in a consistent way.

\begin{lem}
For $a_\nu$ (\ref{4.9})
\label{l5.1}
\begin{equation}
\label{5.1}
(\frac{\partial}{\partial\nu}a_\nu-\frac{\partial}{\partial\nu}a_{-\nu})|_{\nu=n}=
a_n\tilde a,
\end{equation}
where
$$
\tilde a=\frac{2\ln{q^2}}{I_0^{(2)}(2;q^2)}\sum_{k=0}^\infty
\frac{q^{2k^2}}{(q^2,q^2)_k^2}(k-\sum_{l=1}^\infty\frac{q^{2l+2k+2}}{1-q^{2l+2k+2}}).
$$
\end{lem}

{\bf Proof.} Let $n=[\nu]$. It follows from (\ref{4.14})
$$
a_\nu=a_{\nu-n}q^{-n\nu+\frac{n^2}{2}},\qquad
a_{-\nu}=a_{-\nu+n}q^{-n\nu+\frac{n^2}{2}}.
$$
Then

$$
(\frac{\partial}{\partial\nu}a_\nu-
\frac{\partial}{\partial\nu}a_{-\nu})|_{\nu=n}=
(\frac{\partial}{\partial\nu}a_{\nu-n}-
\frac{\partial}{\partial\nu}a_{-\nu+n})|_{\nu=n}
q^{-\frac{n^2}{2}}.
$$
Using (\ref{4.9}) we obtain
$$
(\frac
{\partial}
{\partial\nu}
a_\nu-\frac{\partial}{\partial\nu}
a_{-\nu})|_{\nu=0}=
$$
$$
=\frac
{\sqrt{2/(1-q^2)}e_q(-1)I_0^{(2)}(2;q^2)}
{\Phi_0(2/(1-q^2)}
\frac
{q^{-\frac{n^2}{2}}}
{I_0^{(2)}(2;q^2)}
\left(\frac{\partial}{\partial\nu}I_\nu^{(2)}(2;q^2)
-\frac{\partial}{\partial\nu}I_{-\nu}^{(2)}(2;q^2)\right)|_{\nu=0}=
$$
\begin{equation}
\label{5.3}
=a_0\frac
{q^{-\frac{n^2}{2}}}
{I_0^{(2)}(2;q^2)}
\left(\frac{\partial}{\partial\nu}I_\nu^{(2)}(2;q^2)
-\frac{\partial}{\partial\nu}I_{-\nu}^{(2)}(2;q^2)\right)|_{\nu=0}.
\end{equation}
It follows from (\ref{3.6})

$$
\frac{\partial}{\partial\nu}I_\nu^{(2)}(2;q^2)=\ln{q^2}
\sum_{k=0}^\infty\frac{kq^{2k(\nu+k)}}{(q^2,q^2)_k(1-q^2)^{\nu+k}\Gamma_{q^2}(\nu+k+1)}-
\ln(1-q^2)I_\nu^{(2)}(2;q^2)-
$$
$$
-\sum_{k=0}^\infty
\frac
{q^{2k(\nu+k)}\psi_{q^2}(\nu+k+1)}
{(q^2,q^2)_k(1-q^2)^{\nu+k}\Gamma_{q^2}(\nu+k+1)}.
$$
{}From (\ref{12}) and (\ref{13})
$$
(\frac{\partial}{\partial\nu}I_\nu^{(2)}(2;q^2)
-\frac{\partial}{\partial\nu}I_{-\nu}^{(2)}(2;q^2))|_{\nu=0}=
$$
$$
=2\ln{q^2}\sum_{k=0}^\infty\frac{kq^{2k^2}}{(q^2,q^2)_k^2}-
2\ln(1-q^2)I_0^{(2)}(2;q^2)-
2\sum_{k=0}^\infty\frac{q^{2k^2}}{(q^2,q^2)_k^2}\psi_{q^2}(k+1)=
$$
\begin{equation}
\label{5.4}
=2\ln{q^2}\sum_{k=0}^\infty\frac{kq^{2k^2}}{(q^2,q^2)_k^2}-
2\ln{q^2}\sum_{k=0}^\infty\frac{q^{2k^2}}{(q^2,q^2)_k^2}
\sum_{l=1}^\infty\frac{q^{2l+2k+2}}{1-q^2l+2k+2}.
\end{equation}
Now (\ref{5.1}) follows from (\ref{5.3}) and (\ref{5.4}).

\begin{defi}
\label{d5.1}
The $q$-Bessel-Macdonald functions ($q$-BMF) are defined as
\begin{equation}
\label{5.5}
\makdj=\frac{q^{-\nu^2+1/2}}{4(a_\nu a_{-\nu})^{3/2}\sin{\nu\pi}}
\left[a_\nu\mmbesj-a_{-\nu}\mbesj\right]
\end{equation}
with $a_\nu$ (\ref{4.9}), $j=1, 2$.
\end{defi}

As in the classical case this definition should be adjusted for the integer
values of the index $\nu=n$.
Consider the limit of
(\ref{5.5}) for $j=1$ taking into account (\ref{4.18})
$$
K_n^{(1)}((1-q^2)z;q^2)=\frac{q^{-n^2+1/2}}{4a_n^3}
\lim_{\nu\to n}\frac{a_\nu\mmbes-a_{-\nu}\mbes}{\sin{\nu\pi}}.
$$
Using (\ref{3.5}) we can write
$$
\frac{\partial}{\partial\nu}\mbes=\ln{z/2}\mbes-
\sum_{k=0}^\infty\frac{(1-q^2)^k(z/2)^{\nu+k}\psi_{q^2}(\nu+k+1)}
{(q^2,q^2)_k\Gamma_{q^2}(\nu+k+1)}.
$$
Due to Corollary \ref{c2.1}
$$
[\frac{\partial}{\partial\nu}\mmbes-
\frac{\partial}{\partial\nu}\mbes]|_{\nu=n}=
-2\ln{z/2}I_n^{(1)}((1-q^2)z;q^2)+
$$
$$
+\ln{q^2}\sum_{k=0}^{n-1}(-1)^{n-k-1}
\frac{q^{-(n-k-1)(n-k-2)}(1-q^2)^{-n+2k}(q^2,q^2)_{n-k-1}}
{(q^2,q^2)_k}(z/2)^{-n+2k}+
$$
\begin{equation}
\label{5.6}
+2\sum_{k=0}^\infty\frac{(1-q^2)^{n+2k}(z/2)^{n+2k}}
{(q^2,q^2)_k(q^2,q^2)_{n+k}}
[\psi_{q^2}(n+k+1)+\psi_{q^2}(k+1)].
\end{equation}
Thus
$$
K_n^{(1)}((1-q^2)z;q^2)=\frac{(-1)^nq^{-n^2+1/2}}{4\pi a_n^3}
[(\frac{\partial}{\partial\nu}a_\nu
-\frac{\partial}{\partial\nu}a_{-\nu})|_{\nu=n}I_n^{(1)}((1-q^2)z;q^2)+
$$
$$
+a_n(\frac{\partial}{\partial\nu}\mmbes
-\frac{\partial}{\partial\nu}\mbes)|_{\nu=n}].
$$
The final expression follows from Lemma \ref{l5.1} and (\ref{5.6}) that
$$
K_n^{(1)}((1-q^2)z;q^2)=\frac{(-1)^nq^{-n^2+1/2}}{4\pi a_n^2}
\{(\tilde a-2\ln{z/2})I_n^{(1)}((1-q^2)z;q^2)+
$$
$$
+\ln{q^2}\sum_{k=0}^{n-1}(-1)^{n-k-1}
\frac{q^{-(n-k-1)(n-k-2)}(1-q^2)^{-n+2k}(q^2,q^2)_{n-k-1}}
{(q^2,q^2)_k}(z/2)^{-n+2k}+
$$
$$
+\sum_{k=0}^\infty\frac{(1-q^2)^{n+2k}(z/2)^{n+2k}}
{(q^2,q^2)_k(q^2,q^2)_{n+k}}
\left[\psi_{q^2}(n+k+1)+\psi_{q^2}(k+1)\right]\}.
$$
For $j=2$ we have
$$
K_n^{(2)}((1-q^2)z;q^2)=\frac{(-1)^nq^{-n^2+1/2}}{4\pi a_n^2}
\{(\tilde a-2\ln{z/2})I_n^{(2)}((1-q^2)z;q^2)+
$$
$$
+\ln{q^2}\sum_{k=0}^{n-1}(-1)^{n-k-1}
\frac{q^{-(n-k-1)(n-k-2)+2k(k-n)}(1-q^2)^{-n+2k}(q^2,q^2)_{n-k-1}}
{(q^2,q^2)_k}(z/2)^{-n+2k}+
$$
$$
+\sum_{k=0}^\infty\frac{q^{2k(n+k)}(1-q^2)^{n+2k}(z/2)^{n+2k}}
{(q^2,q^2)_k(q^2,q^2)_{n+k}}
\left[\psi_{q^2}(n+k+1)+\psi_{q^2}(k+1)\right]\}.
$$

\begin{predl}
\label{p5.1}
q-BMF $\makd$ is a holomorphic
function in the domain  $Re z>\frac{2q}{1-q^2}$
\end{predl}

{\bf Proof.} Substitute (\ref{4.8}) in (\ref{5.5}). Then
\begin{equation}
\label{5.8}
\makd=\frac{q^{-\nu^2+1/2}}{2\sqrt{a_\nu a_{-\nu}}\sqrt{z}}
\mqexp\Phi_\nu(z).
\end{equation}
The product of
$\mqexp$ and
$\basgyp$ is a holomorphic function in Re$z>\frac{2q}{1-q^2}$.

\begin{predl}
\label{p5.1a}
q-BMF $\makdo$ is a holomorphic
function in the domain  $z\ne0$
\end{predl}

{\bf Proof.} This statement follows from (\ref{5.5}) and (\ref{3.6}).

\begin{cor}
\label{c5.1}
\begin{equation}
\label{5.8a}
\makd0=\frac{q^{-\nu^2+1/2}}{2\sqrt{a_\nu a_{-\nu}}\sqrt{z}}
E_q(-\frac{1-q^2}{2}z)\Phi_\nu(z).
\end{equation}
\end{cor}

{\bf Proof.} This formula follows from (\ref{5.5}) and (\ref{4.19}).

\begin{predl}
\label{p5.2}
The function $\makd$ satisfies the following
relations
\begin{equation}
\label{5.9}
\frac{2}{(1+q)z}\partial_qz^\nu\makd=-z^{\nu-1}K_{\nu-1}((1-q^2)z;q^2),
\end{equation}
$$
\frac{2}{(1+q)z}\partial_qz^{-\nu}\makd
=-z^{-\nu-1}K_{\nu+1}((1-q^2)z;q^2).
$$
\end{predl}

{\bf Proof.} It follows from Definition 5.1 and Proposition \ref{p3.1}
$$
\frac{2}{(1+q)z}\partial_qz^\nu\makd=
\frac{q^{-\nu^2+1/2}}{4(a_\nu a_{-\nu})^{3/2}\sin{\nu\pi}}
\left[a_\nu\frac{2}{(1+q)z}\partial_qz^\nu I_{-\nu}^{(1)}-
a_{-\nu}\frac{2}{(1+q)z}\partial_qz^\nu I_\nu^{(1)}\right]=
$$
$$
=\frac{q^{-\nu^2+1/2}}{4(a_\nu a_{-\nu})^{3/2}\sin{\nu\pi}}
[a_\nu z^{\nu-1}I_{-\nu+1}^{(1)}-a_{-\nu}z^{\nu-1}I_{\nu-1}^{(1)}].
$$
As it follows from (\ref{4.14}) $a_\nu=a_{\nu-1}q^{-\nu+1/2},\qquad
a_{-\nu}=a_{-\nu+1}q^{-\nu+1/2}$. So
$$
\frac{2}{(1+q)z}\partial_qz^\nu K_\nu^{(1)}=
$$
$$
=-z^{\nu-1}
\frac{q^{-(\nu-1)^2+1/2}}{4(a_{\nu-1}a_{-\nu+1})^{3/2}\sin{(\nu-1)\pi}}
[a_{\nu-1}
z^{\nu-1}I_{-\nu+1}^{(1)}-a_{-\nu+1}z^{\nu-1}I_{\nu-1}^{(1)}]=
-z^{\nu-1}K_{\nu-1}^{(1)}.
$$
It is easy to see that
$K_{-\nu}^{(j)}((1-q^2)z;q^2)=K_\nu^{(j)}((1-q^2)z;q^2)$.
Thus from (\ref{5.9})
$$
\frac{2}{(1+q)z}\partial_qz^{-\nu}K_\nu^{(1)}=
\frac{2}{(1+q)z}\partial_qz^{-\nu}K_{-\nu}^{(1)}=
-z^{-\nu-1}K_{-\nu-1}^{(1)}=
-z^{-\nu-1}K_{\nu+1}^{(1)}.\rule{5pt}{5pt}
$$

\begin{predl}
\label{p5.3}
The function $\makd$ satisfies the next functional
relations
$$
K_{\nu-1}^{(1)}((1-q^2)z;q^2)-K_{\nu+1}^{(1)}((1-q^2)z;q^2)=
-\frac{2}{(1-q^2)z}(q^{-\nu}-q^\nu)K_\nu^{(1)}((1-q^2)qz;q^2),
$$
$$
K_{\nu-1}^{(1)}((1-q^2)z;q^2)+K_{\nu+1}^{(1)}((1-q^2)z;q^2)=
-\frac{4}{(1-q^2)z}\makd+
$$
$$
+\frac{2}{(1-q^2)z}(q^{-\nu}+q^\nu)K_\nu^{(1)}((1-q^2)qz;q^2).
$$
\end{predl}

{\bf Proof.} Denote the coefficients in front of the square brackets in
(\ref{5.5}) by $A_\nu$. Then  from (\ref{4.14})
$$
A_{\nu-1}=-A_\nu q^{-\nu+1/2},\qquad A_{\nu+1}=-A_\nu q^{\nu+1/2}.
$$
The statement  follows from
(\ref{4.14}),(\ref{5.5}),(\ref{3.9}) and (\ref{3.10}).

\begin{predl}
\label{p5.2a}
The function $\makdo$ satisfies the following
relations
$$
\frac{2}{(1+q)z}\partial_qz^\nu\makdo=
-q^{-\nu+1}z^{\nu-1}K_{\nu-1}^{(2)}((1-q^2)z;q^2),
$$
$$
\frac{2}{(1+q)z}\partial_qz^{-\nu}\makdo
=-q^{-\nu+1}z^{-\nu-1}K_{\nu+1}^{(2)}((1-q^2)z;q^2).
$$
\end{predl}

\begin{predl}
\label{p5.3a}
The function $\makdo$ satisfies the next functional
relations
$$
K_{\nu-1}^{(2)}((1-q^2)z;q^2)-K_{\nu+1}^{(2)}((1-q^2)z;q^2)=
-\frac{2q^{\nu-1}}{(1-q^2)z}(q^{-\nu}-q^\nu)K_\nu^{(2)}((1-q^2)qz;q^2),
$$
$$
K_{\nu-1}^{(2)}((1-q^2)z;q^2)+K_{\nu+1}^{(2)}((1-q^2)z;q^2)=
-\frac{4q^{\nu-1}}{(1-q^2)z}\makdo+
$$
$$
+\frac{2q^{\nu-1}}{(1-q^2)z}(q^{-\nu}+q^\nu)K_\nu^{(2)}((1-q^2)qz;q^2).
$$
\end{predl}

The proof of Propositions \ref{p5.2a} and \ref{p5.3a} are the same as
\ref{p5.2} and \ref{p5.3}.

\begin{rem}
\label{r5.1}
\begin{equation}
\label{5.13}
\lim_{q\to 1-0}\makdj=K_\nu(z),\qquad j=1, 2.
\end{equation}
\end{rem}

Really it follows from (\ref{4.14}) and (\ref{4.15}) if $q=1$ then
$a_\nu$ is independent of $\nu$, and $a_\nu=\frac{1}{\sqrt{2\pi}}$. Now
(\ref{5.13}) follows from Remark \ref{r3.1}.

\begin{rem}
\label{r5.2}
If $q\to 1-0$ the representations (\ref{4.8}) and
(\ref{5.8}) give us the well known asymptotic decompositions for the
functions $I_\nu(z)$ and $K_\nu(z)$ respectively \cite{BE}.
\end{rem}

\section{Some relations for the basic hypergeometric functions}
Here we summarize some relations for the basic hypergeometric functions that
we have already derived in Section 4. It is worthwhile to note that they are
based
on the identifcations of the power and Laurent expansions for the same
$q$-functions. In the limit $q\rightarrow 1$ it is became impossible since the
later
expansion  diverges. Remind that
\begin{equation}
\label{6.1}
\basgyp=\Phi_\nu(z),
\end{equation}
$$
a_\nu=\sqrt{\frac{2}{1-q^2}}e_q(-1)
\frac{I_\nu^{(2)}(2;q^2)}{\Phi_\nu(\frac{2}{1-q^2})},~\nu\neq n,~~
a_n=\sqrt{\frac{q^{-n^2+1/2}\ln{q^{-2}}}{2\pi(1-q^2)}}.
$$
Assume in(\ref{6.1}) $\frac{2q}{(1-q^2)z}=u$,\quad$|u|<1$.
Then
$$
a_{\nu-1}\phantom._2\Phi_1(q^{\nu-1/2},q^{-\nu+3/2};-q;q,u)
-a_{\nu+1}\phantom._2\Phi_1(q^{\nu+3/2},q^{-\nu-1/2};-q;q,u)=
$$
$$
=a_{\nu}q^{-1/2}(q^{-\nu}-q^\nu)\left(u/q-1\right)
\phantom._2\Phi_1(q^{\nu+1/2},q^{-\nu+1/2};-q;q,u/q),
$$
$$
a_{\nu-1}\phantom._2\Phi_1(q^{\nu-1/2},q^{-\nu+3/2};-q;q,u)
+a_{\nu+1}\phantom._2\Phi_1(q^{\nu+3/2},q^{-\nu-1/2};-q;q,u)=
$$
$$
=2a_\nu u/q\phantom._2\Phi_1(q^{\nu+1/2},q^{-\nu+1/2};-q;q,u)-
$$
$$
-a_{\nu}q^{-1/2}(q^{-\nu}+q^\nu)\left(u/q-1\right)
\phantom._2\Phi_1(q^{\nu+1/2},q^{-\nu+1/2};-q;q,u/q).
$$

The last relation coming from the q-Wronskian (\ref{3.17}) has been proved for
the nonineteger $\nu$
$$
(u/q+1)\phantom._2\Phi_1(q^{\nu+1/2},q^{-\nu+1/2};-q;q,u)
\phantom._2\Phi_1(q^{\nu+1/2},q^{-\nu+1/2};-q;q,-u/q)-
$$
$$
(u/q-1)\phantom._2\Phi_1(q^{\nu+1/2},q^{-\nu+1/2};-q;q,-u)
\phantom._2\Phi_1(q^{\nu+1/2},q^{-\nu+1/2};-q;q,u/q)=2.
$$

Since $\phantom._2\Phi_1$ is the continuous function of $\nu$, then
this equality  is valid for integer $\nu$.

\small{

}

\end{document}